\documentclass[11pt,amsfonts]{article}
\usepackage{graphicx,amssymb,eucal,mathrsfs}
\usepackage{latexsym,amsmath,amscd,epsfig,xcolor}
\usepackage[all]{xy}
\usepackage{ntheorem,tikz}
\usepackage{epsf,graphicx,pgf,etex,amscd,pstricks}

\usepackage[hypertexnames=false,backref=page,pdftex,
 	pdfpagemode=UseNone,
 	breaklinks=true,
 	extension=pdf,
 	colorlinks=true,
 	linkcolor=blue,
 	urlcolor=blue,
 ]{hyperref}

\newcommand{\Z}{\mathbb{Z}}
\newcommand{\Q}{\mathbb{Q}}
\newcommand{\R}{\mathbb{R}}

\newcommand{\Hy}{\mathbb{H}}
\newcommand{\Xs}{\mathbb{X}}

\newcommand{\N}{\mathbb{N}}
\newcommand{\del}{{\partial}}

\newcommand{\PO}{{\rm PO}(n,1)}

\newcommand{\lk}{{\rm lk}}

\newtheorem{main}{Theorem}

\newtheorem{question}[main]{Question}
\newtheorem{prop}[main]{Proposition}
\newtheorem{lemma}[main]{Lemma}
\newtheorem{cor}[main]{Corollary}

\theorembodyfont{\upshape}

\newtheorem{defi}[main]{Definition}
\newtheorem{rem}[main]{Remark}
\newtheorem{exam}[main]{Example}

\renewcommand{\thefootnote}{\alph{footnote}}

\title{Hyperbolic groups containing subgroups of type $\mathscr{F}_{3}$ not $\mathscr{F}_{4}$}
\date{2021}
\author{Claudio Llosa Isenrich$^{{\rm a}}$, Bruno Martelli and Pierre Py}

\begin{document}

\maketitle

\begin{abstract} We give examples of hyperbolic groups which contain subgroups that are of type $\mathscr{F}_{3}$ but not of type $\mathscr{F}_{4}$. These groups are obtained by Dehn filling starting from a non-uniform lattice in ${\rm PO}(8,1)$ which was previously studied by Italiano, Martelli and Migliorini. 
\end{abstract}

\footnotetext[1]{The author gratefully acknowledges funding by the DFG 281869850 (RTG 2229)}

\renewcommand{\thefootnote}{\arabic{footnote}}


\section{Introduction}
\label{sec:Intro}

The purpose of this work is to prove the following theorem. 

\begin{main}\label{ssgrofhyperbolicgroup} There exist infinitely many pairwise non-isomorphic hyperbolic groups $G$ admitting a surjective homomorphism $\phi \colon G \to \Z$ whose kernel is of type $\mathscr{F}_{3}$ and has the property that $H_{4}(\ker(\phi),\Z)$ is not finitely generated. In particular $\ker(\phi)$ is not of type $\mathscr{F}_{4}$. 
\end{main}

These hyperbolic groups will be obtained by {\it Dehn filling} starting from a certain non-uniform lattice in ${\rm PO}(8,1)$. Before describing our method, we provide some historical context for our result.

The study of finiteness properties of subgroups of hyperbolic groups has a long history. Rips gave the first example of a finitely generated subgroup of a hyperbolic group which is not finitely presented~\cite{rips}. This was followed by many works on {\it coherence}, see e.g.~\cite{bome,JNW,kapovich,kapovichpotyagailo,kapovi} to name just a few of them. While these works lead to an ample supply of finitely generated subgroups of hyperbolic groups that are not themselves hyperbolic, finding finitely presented subgroups of hyperbolic groups that are not themselves hyperbolic turns out to be a much more subtle problem. 

The first example of a subgroup of a hyperbolic group which is finitely presented but not of type $\mathscr{F}_{3}$, and thus not itself hyperbolic, was constructed by Brady~\cite{brady}  (see also~\cite{bradyrileyshort} for a survey). Preceding Brady's work, Gromov~\cite{gromov} was the first to sketch the construction of a ramified covering of a 5-dimensional torus that seemed to produce a non-hyperbolic finitely presented subgroup of a hyperbolic group. However, it was later shown by Bestvina that Gromov's group cannot be hyperbolic \cite{bestvina}. Brady's construction in \cite{brady} elegantly solves this problem, by taking a covering of a direct product of three graphs, ramified on a family of embedded graphs. More recently, other examples with the same property as Brady's group were constructed by Lodha~\cite{lodha} and by Kropholler~\cite{krophollergardam}, who generalises both Brady's and Lodha's constructions. Finally, the very first examples of subgroups of hyperbolic groups which are not hyperbolic but have finite classifying spaces were constructed by Italiano, Martelli and Migliorini in 2021~\cite{itmami2}. The existence of such subgroups was a well-known open problem before their work. A variation of their construction was subsequently provided by \cite{fujiwara}. To the best of our knowledge this provides a complete overview on known examples of non-hyperbolic finitely presented subgroups of hyperbolic groups.

In \cite{brady}, Brady raised the following question. 

\begin{question}
Do there exist hyperbolic groups containing subgroups of type $\mathscr{F}_{n-1}$ but not of type $\mathscr{F}_n$ for $n\ge 4$?
\end{question}

Theorem~\ref{ssgrofhyperbolicgroup} gives a positive answer to Brady's question for $n=4$. We emphasize that parts of the methods we develop work in all dimensions and could thus also lead to examples of subgroups of hyperbolic groups of type $\mathscr{F}_{n-1}$ and not of type $\mathscr{F}_{n}$ for $n\geq 5$. This would require finding examples of non-uniform lattices in ${\rm PO}(2n,1)$ endowed with homomorphisms onto $\Z$ whose kernels have strong enough finiteness properties. 

Some evidence towards the existence of such lattices is provided by a recent theorem by Fisher~\cite{fisher}, proving a conjecture of Kielak~\cite{kielakoberwolfach}. Kielak's conjecture was motivated by~\cite{kielak}, and one of the consequences of its proof is the existence of subgroups of hyperbolic groups which are of type ${\rm FP}_{n-1}(\mathbb{Q})$ but not of type ${\rm FP}_{n}(\mathbb{Q})$ for each $n\ge 3$.  We discuss this briefly in Section~\ref{nonfkness}. Even though this does not seem to provide new examples of subgroups of hyperbolic groups of type $\mathscr{F}_{n-1}$ and not $\mathscr{F}_{n}$ for $n\geq 3$, it does suggest in conjunction with Theorem \ref{ssgrofhyperbolicgroup} that lattices in ${\rm PO}(2n,1)$ are a good place to search for such groups.

We now describe the construction leading to the examples of Theorem~\ref{ssgrofhyperbolicgroup}. Let $\Gamma < \PO$ be a torsion-free nonuniform lattice where we assume that $n\ge 3$. Let $(H_{i})_{i\in I}$ be the family of maximal parabolic subgroups of $\Gamma$. We assume that each $H_i$ is purely unipotent, hence is isomorphic to $\Z^{n-1}$. We fix once and for all a collection $(B_{i})_{i\in I}$ of open horoballs in the $n$-dimensional hyperbolic space $\mathbb{H}^{n}$, in such a way that: 
\begin{enumerate}
\item $B_i$ is invariant under $H_i$ for each $i\in I$,
\item the $B_i$'s are disjoint and the union $\cup_{i\in I}B_{i}$ is $\Gamma$-invariant,
\item the action of $\Gamma$ on the closed set $\widetilde{C}=\mathbb{H}^{n}-\cup_{i\in I}B_{i}$ is cocompact.  
\end{enumerate}
If $\Lambda$ is a finite index subgroup of $\Gamma$ we define $M_{\Lambda}:=\widetilde{C}/\Lambda$. Hence
\[
M_{\Lambda}\subset \mathbb{H}^{n}/\Lambda
\]
is a compact submanifold whose boundary is a disjoint union of finitely many flat $(n-1)$-tori obtained by cutting the cusps of $\mathbb{H}^{n}/\Lambda$ along the projection of the horoballs $B_i$.

Let now $f \colon \Gamma \to \Z$ be a surjective homomorphism. Throughout the introduction, we make the following:

\noindent {\bf Standing assumption.} The homomorphism $f$ is nontrivial on each parabolic subgroup. 

If $\Lambda < \Gamma$ is a finite index subgroup, we denote by $N_{f}(\Lambda)$ the normal subgroup generated by the union:
\[
\bigcup_{i\in I} H_{i}\cap \Lambda \cap \ker (f)
\]
and by $\overline{f}_{\Lambda} \colon \Lambda /N_{f}(\Lambda)\to \Z$ the homomorphism induced by $f$. The groups $H_{i}\cap \Lambda$ are isomorphic to $\Z^{n-1}$ and each intersection
\[
H_{i}\cap \Lambda \cap \ker(f)
\]
is isomorphic to $\Z^{n-2}$. Geometrically, $\Lambda/N_{f}(\Lambda)$ is the fundamental group of the space $M_{\Lambda,f}$ obtained from $M_{\Lambda}$ by coning off in each boundary component a family of codimension $1$ subtori parametrised by $S^1$. These subtori are determined by the condition that their fundamental group coincides with the kernel of the restriction of $f$ to the fundamental group of the corresponding boundary component.

We will deduce from a result due to Fujiwara and Manning~\cite{fujiwaramanning} about Dehn filling in cusped hyperbolic manifolds that the quotient $\Lambda /N_{f}(\Lambda)$ is hyperbolic if $\Lambda$ is deep enough in $\Gamma$. We will describe this in more detail and give a precise definition of the expression ``deep enough" in Section~\ref{dehnfi}. The hyperbolicity of $\Lambda /N_{f}(\Lambda)$ can also be deduced from the more general theory of Dehn filling in relatively hyperbolic groups~\cite{dgo,grovesmanning,osin}. Our first result is the following:  

\begin{main}\label{fkispreserved} Let $k\ge 1$. If the kernel of $f$ is of type $\mathscr{F}_{k}$, then the kernel of the induced morphism $\overline{f}_{\Lambda} \colon \Lambda/N_{f}(\Lambda)\to \Z$ is also of type $\mathscr{F}_{k}$ for every deep enough finite index subgroup $\Lambda < \Gamma$. 
\end{main}

The main input to prove Theorem~\ref{fkispreserved} will be Fujiwara and Manning's result~\cite[Theorem 2.7]{fujiwaramanning} stating that the coned-off space $M_{\Lambda,f}$ carries a locally ${\rm CAT}(-1)$ metric (for deep enough $\Lambda$), which implies that $M_{\Lambda,f}$ is aspherical. We spell out the particular case of Theorem~\ref{fkispreserved} when $k=2$:

\begin{main}\label{finitepresentationispreserved} If the kernel of $f$ is finitely presented, then the kernel of the induced morphism $\overline{f}_{\Lambda} \colon \Lambda/N_{f}(\Lambda)\to \Z$ is also finitely presented for every deep enough finite index subgroup $\Lambda < \Gamma$.  
\end{main}

Theorem \ref{fkispreserved} leads to two natural questions. Can one find examples of pairs $(\Gamma, f)$ to which one can apply this construction? Given such a pair, can one determine the exact finiteness properties of the kernel $\ker (\overline{f}_{\Lambda})$ associated to a deep enough finite index subgroup $\Lambda < \Gamma$? We first give a partial answer to the second question when $n$ is even.

\begin{main}\label{theo:nonfiniteness} Assume that $n=2k$ and that the Betti numbers $(b_{i}(\ker (f))_{1\le i \le k-1}$ are finite. Then $b_{k}(\ker(f))$ is infinite. If, moreover, $\Lambda$ is deep enough, then $b_{k}(\ker(\overline{f}_{\Lambda}))$ is also infinite and in particular $\ker(\overline{f}_{\Lambda})$ is not of type $\mathscr{F}_k$.
\end{main} 

The proof of Theorem~\ref{theo:nonfiniteness} relies on classical arguments due to Milnor~\cite{milnor}, combined with the fact that the Euler characteristic of finite volume (complete, oriented) even-dimensional hyperbolic manifolds is nonzero~\cite{kz}. 

To deduce Theorem~\ref{ssgrofhyperbolicgroup} from Theorems~\ref{fkispreserved} and~\ref{theo:nonfiniteness}, we will use a certain hyperbolic manifold $M^8=\mathbb{H}^{8}/\Gamma$ built from the dual of the 8-dimensional Euclidean Gosset polytope studied in~\cite{everittratcliffetschantz}, together with a continuous map $u \colon M^{8} \to S^{1}$. The manifold $M^{8}$ and the corresponding map to the circle were constructed in~\cite{itmami}. There the authors prove that the kernel of the homomorphism 
$$u_{\ast} \colon \Gamma \to \Z$$
induced by $u$ is finitely presented. Here we go further and prove that $\ker (u_{\ast})$ is actually of type $\mathscr{F}_{3}$ and not of type $\mathscr{F}_{4}$ (see Theorem~\ref{thm:Finprops-M8}). However, the restriction of the map $u$ to some of the cusps is homotopic to a constant function. By considering a suitable rational perturbation of the cohomology class of $u$, we prove:

\begin{main}\label{goodmaptothecircle} There exists a continuous map $v \colon M^{8} \to S^{1}$ such that:
\begin{enumerate}
\item none of the restrictions of $v$ to a cusp of $M^{8}$ is homotopic to a constant
\item the kernel of $v_{\ast}$ is of type $\mathscr{F}_{3}$.
\end{enumerate} 
\end{main}

For the part of our work concerning the specific properties of this 8-dimensional manifold we use a code written in Sage to perform some of the computations required in our argument.

Combining Theorems~\ref{fkispreserved}, \ref{theo:nonfiniteness} and~\ref{goodmaptothecircle}, and applying them to $f=v_{\ast}$, we can prove the existence of one hyperbolic group satisfying the conclusion of Theorem~\ref{ssgrofhyperbolicgroup}: if $\Lambda<\Gamma$ is a deep enough finite index subgroup, the fundamental group $\Lambda/N_{f}(\Lambda)$ of the coned off space $M_{\Lambda,f}$ is hyperbolic and the kernel of the induced morphism
$$\overline{f}_{\Lambda} \colon \Lambda/N_{f}(\Lambda) \to \Z$$
has the properties stated in Theorem~\ref{ssgrofhyperbolicgroup}.

The only point that remains to be proved is then the fact that one can obtain from our methods infinitely many isomorphism classes of hyperbolic groups containing subgroups of type $\mathscr{F}_3$ and not $\mathscr{F}_4$. This will follow from a result due to Fujiwara and Manning~\cite[Theorem 2.13]{fujiwaramanning} by constructing a suitable infinite family of perturbations of $u_{\ast}$ as in Theorem \ref{goodmaptothecircle}. This will require a few additional arguments and we will postpone this part of the proof to Section \ref{sec:infinitefamily}.


\section{Dehn fillings and finiteness properties}\label{dehnfillingfiniteness}

In this section we will prove Theorem~\ref{fkispreserved}. For this we first recall Fujiwara and Manning's results on Dehn fillings~\cite{fujiwaramanning} of non-compact finite volume hyperbolic manifolds in Section~\ref{dehnfi}. In Section ~\ref{section:finiteness} we then introduce different types of finiteness properties of groups and prove Theorem~\ref{fkispreserved}.

\subsection{Dehn filling and asphericity}\label{dehnfi}

We keep the notations from the introduction: $\Gamma$ is a torsion-free nonuniform lattice of isometries of $\mathbb{H}^{n}$, the family of its parabolic subgroups is denoted by $(H_{i})_{i\in I}$ and we have fixed a family $(B_{i})_{i\in I}$ of horoballs having the same properties as in the introduction. The complement of the union of the $B_{i}$'s in $\mathbb{H}^{n}$ is denoted by $\widetilde{C}$ and we set $M_{\Gamma}=\widetilde{C}/\Gamma$. We pick a finite set
$$\{i_{1}, \ldots , i_{r}\}\subset I$$
with the property that each parabolic subgroup of $\Gamma$ is in the orbit of exactly one of the $(H_{i_{s}})_{1\le s\le r}$. The boundary of $B_{i_{s}}$ projects onto a boundary torus $T_{s}$ in $M_{\Gamma}$. The metric induced by the hyperbolic metric on $T_{s}$ is flat. We now describe more precisely the coned-off space $M_{\Gamma,f}$ (see~\cite[\S 2.1]{fujiwaramanning} for a similar, more general description). Identify $T_s$ with its flat metric with $\R^{n-1}/A_{s}$ for a lattice $A_{s}<\R^{n-1}$. Let $D_{s}<A_{s}$ be the intersection of the kernel of $f$ with $A_{s}$. The span of $D_s$ is a codimension $1$ subspace $V\subset \R^{n-1}$ and the orthogonal projection $\R^{n}\to V^{\perp}$ induces a submersion $\pi_s \colon T_{s}\to S^{1}$ whose fibers are the translates of the subtorus $V/D_{s}$. The {\it partial cone} $C(T_{s},V/D_{s})$ is the space obtained from $T_{s}\times [0,1]$ by identifying $(p,1)$ to $(q,1)$ whenever $\pi_s(p)=\pi_s (q)$. The map $(p,t)\mapsto \pi_{s}(p)$ induces a map
\begin{equation}\label{conetocircle}
C(T_{s},V/D_{s})\to S^{1}.
\end{equation}
The space $M_{\Gamma,f}$ is obtained by gluing $C(T_{s},V/D_{s})$ to $M_{\Gamma}$, identifying $T_{s}\times \{0\}$ with $T_{s}$ for each $s\in \left\{1,\dots,r\right\}$. Clearly, one can perform this construction replacing $\Gamma$ by any of its finite index subgroups.

We say that a torus equipped with a flat metric satisfies the $2\pi$-condition if each of its closed geodesics has length greater than $2\pi$. Similarly, we say that a totally geodesic subtorus of a flat torus satisfies the $2\pi$-condition if each of the closed geodesics in the subtorus has length greater than $2\pi$. Clearly if a torus satisfies the $2\pi$-condition all of its totally geodesic subtori also satisfy it. The following theorem follows from Theorem 2.7 in~\cite{fujiwaramanning}. 

\begin{main}\label{theo:fujiwaramanning} Assume that for each $s\in \{1, \ldots , r\}$, $V/D_s$ satisfies the $2\pi$-condition. Then $M_{\Gamma, f}$ carries a locally ${\rm CAT}(-1)$ metric.   
\end{main}

A priori the $2\pi$-condition for $V/D_s$ need not be satisfied. However, there is always some finite index subgroup $\Lambda<\Gamma$ such that the boundary components of $M_{\Lambda}$ satisfy the $2\pi$-condition. We recall this classical argument. For $1\le s \le r$, let $Z_{s}\subset H_{i_{s}}$ be the set of elements whose translation length for the action $H_{i_{s}} \curvearrowright \partial B_{i_{s}}$ is not greater than $2\pi$. This is a finite set. Hence
$$Z\colon=Z_{1}\cup \cdots \cup Z_{r}\subset \Gamma$$
is a finite set.

\begin{defi} A finite index subgroup $\Lambda < \Gamma$ is \emph{deep enough} if the intersection of $\Lambda$ with the union $\cup_{\gamma \in \Gamma}\gamma Z\gamma^{-1}$ is trivial.  
\end{defi}
 
 Deep enough finite index subgroups exist: indeed, since $\Gamma$ is residually finite, there is a finite quotient $a \colon \Gamma \to F$ such that $a(z)\neq 1$ for all $z\in Z$. The kernel of such a morphism $a$ is a deep enough finite index subgroup. A possibly nonnormal finite index subgroup of $\ker (a)$ also does the job. As a consequence of Theorem~\ref{theo:fujiwaramanning}, applied to a finite index subgroup of $\Gamma$ instead of $\Gamma$ itself, we have: 

\begin{cor}\label{corollaryoffujiwaramanning} If $\Lambda$ is a deep enough finite index subgroup of $\Gamma$, then the space $M_{\Lambda,f}$ is aspherical and the group $\Lambda/N_{f}(\Lambda)$ is hyperbolic.  
\end{cor}

 Hence, if $\Lambda$ is a deep enough finite index subgroup of $\Gamma$, each covering space $Y\to M_{\Lambda,f}$ is a classifying space for the corresponding subgroup of $\Lambda/N_{f}(\Lambda)$. This will allow us to study the finiteness properties of $\pi_{1}(Y)$ geometrically.

\begin{rem} The asphericity of $M_{\Lambda, f}$ could also be established using the methods from~\cite{coulon}. That paper establishes the asphericity of certain spaces obtained by coning off totally geodesic submanifolds in closed negatively curved manifolds and the same techniques can also be applied in the context of Dehn fillings. 
\end{rem}

\subsection{Finiteness properties}\label{section:finiteness}

In this Section, we recall the definitions of all the finiteness properties that we shall deal with and then prove Theorem~\ref{fkispreserved}. We start with the definition of property $\mathscr{F}_{n}$, introduced by Wall~\cite{wall}.

\begin{defi}
A group $G$ is of type $\mathscr{F}_{n}$ if it admits a $K(G,1)$ which is a CW-complex with finite $n$-skeleton. 
\end{defi}

The above definition is a {\it geometric} finiteness property. We now define some {\it algebraic} finiteness properties, first introduced by Bieri~\cite{bieri}. 

\begin{defi}
Let $R$ be a ring. A group $G$ is of type ${\rm FP}_{n}(R)$ if $R$, considered as a trivial $RG$-module, admits a projective resolution 
$$\cdots \to P_{n+1}\to P_{n} \to \cdots \to P_{1} \to P_{0} \to R \to 0$$
where $P_i$ is a finitely generated $RG$-module for $i\le n$. 
\end{defi}

\noindent In what follows, we will only consider the case when $R=\Z$ (in which case we abbreviate and write that $G$ is of type ${\rm FP}_{n}$ instead of ${\rm FP}_{n}(\Z)$) and the case $R=\Q$. If $G$ is of type ${\rm FP}_{n}(R)$, one can actually construct a resolution 
$$\cdots \to F_{n+1}\to F_{n} \to \cdots \to F_{1} \to F_{0} \to R \to 0$$
by {\it free} $RG$-modules $F_i$, where $F_i$ is finitely generated for $i\le n$~\cite[Proposition VIII.4.3]{brown}. One has the implications
$$\mathscr{F}_{n} \Longrightarrow {\rm FP}_{n} \Longrightarrow {\rm FP}_{n}(\Q).$$
Furthermore, if a group $G$ is of type ${\rm FP}_{n}(\Q)$, its Betti numbers $b_{i}(G)$ are finite for $i\le n$, where $b_{k}(G)$ is the $k$-th Betti number of $G$ with coefficients in $\Q$. We also mention that all homology groups that we will consider will be taken with rational coefficients, if not specified otherwise. 

We now return to the proof of Theorem~\ref{fkispreserved}. We assume that $\ker(f)$ is $\mathscr{F}_k$ and that we fixed a deep enough finite index subgroup $\Lambda<\Gamma$. In particular, the manifold $M_{\Lambda}$, together with the subtori of its boundary determined by ${\rm ker}(f)$, satisfies the $2\pi$-condition. To simplify the presentation, we will use the following notations: $W:=M_{\Lambda}$ is the hyperbolic manifold with boundary, $v\colon W\to S^1$ is a continuous map inducing the morphism $f \colon \Lambda \to \Z$ (obtained by restriction of $f \colon \Gamma \to \Z$), and $W_f$ is the Dehn filling of $W$ constructed in Section \ref{dehnfi}. The induced map and homomorphism are denoted by $\overline{v}\colon W_f\to S^1$ and $\overline{f}\colon\pi_1(W_f)\to \Z$. After passing to the finite covers of $W$ and $W_f$ corresponding to the preimage under $f$ (or $\overline{f}$) of the finite index subgroup $$\bigcap_{1\leq s \leq r} f(A_s)\leq \Z$$ if necessary, we may moreover assume that $f|_{A_s}\colon A_s\to \pi_1(S^1)\cong \Z$ is surjective for $1\leq s \leq r$. Since $\ker (f)$ and $\ker (\overline{f})$ remain the same under these finite covers of $W$ and $W_f$, it suffices to consider this case. Finally we can and do assume that the map $\overline{v}$ coincides on $C(T_{s},V/D_{s})$ with the natural map $C(T_{s},V/D_{s})\to S^{1}$ induced by $\pi_{s}$, see~\eqref{conetocircle}. (This is possible after performing a suitable homotopy of the original map $v$.) 

Let $\pi \colon W_{\Z}\to W$ and $\overline{\pi} \colon W_{f,\Z}\to W_{f}$ be the covering spaces induced by $\ker f$ and $\ker \overline{f}$. Our assumptions imply that  for $1\leq s \leq r$ the preimage $\overline{\pi}^{-1}(C(T_s,V/D_s))$ is connected and can canonically be identified with $C(V/D_s)\times \R$ via the $\R$-translates of the inclusion of the cone $C(V/D_s)\stackrel{\cong}{\rightarrow} v^{-1}(0)\cap C(T_s,V/D_s)$ of $V/D_s$ in $W_{f,\Z}$. We obtain a commutative diagram
\[
\xymatrix{ 
W_{\Z} \ar@{_{(}->}[d] \ar[r]& \R \ar[r]\ar[d]^{{\rm id}_{\R}} &S^1\ar[d]^{{\rm id}_{S^1}}\\
W_{f,\Z}  \ar[r]& \R \ar[r] & S^1\\
\bigcup_{1\leq s \leq r} C(V/D_s)\times \R \ar@{^{(}->}[u]\ar[r] & \R \ar[r]\ar[u]_{{\rm id}_{\R}} & S^1\ar[u]_{{\rm id}_{S^1}},
}
\]
with $W_{f,\Z}\setminus \left(\bigcup_{1\leq s \leq r} C(V/D_s)\times \R\right) = W_{\Z}\setminus \partial W_{\Z}$.  It follows that there is a deformation retraction of the aspherical space $W_{f,\Z}$ onto the subspace $$W_{\Z}\cup \left( \bigcup_{1\leq s \leq r} C(V/D_s)\times \left\{0\right\}\right).$$ In particular, up to homotopy equivalence, $W_{f,\Z}$ is an aspherical space that is obtained by attaching finitely many cells to the boundary of $W_{\Z}$. On the other hand, since the aspherical manifold $W_{\Z}$ is a $K(\ker (f),1)$ for the group $\ker(f)$ of type $\mathscr{F}_k$, it is homotopy equivalent to an aspherical CW-complex with finitely many cells of dimension $\leq k$.
We deduce that $W_{f,\Z}$ is homotopy equivalent to an aspherical CW-complex with finitely many cells of dimension $\leq k$. Thus, $\pi_1(W_{f,\Z})=\ker(\overline{f})$ is of type $\mathscr{F}_k$. This completes the proof of Theorem \ref{fkispreserved}.

We end this section by proving the following variation of Theorem \ref{fkispreserved}, which we will require later.
\begin{prop}\label{prop:Betti-finiteness}
 Let $i\ge 0$. Let $\Lambda < \Gamma$ be a deep enough finite index subgroup, $f\colon\Lambda \to \Z$ and $\overline{f}_{\Lambda}\colon\Lambda/N_f(\Lambda)\to \Z$ be as in Section \ref{dehnfi}. Then $b_i(\ker(f))$ is finite if and only if $b_i({\ker(\overline{f}_{\Lambda}}))$ is finite.
\end{prop}
{\noindent \textit{Proof.}} We make the same assumptions and use the same notation as in the remainder of this section. From the decomposition
\[
W_{f,\Z}= W_{\Z}\bigcup \left(\bigcup_{1\leq s \leq r} C(V/D_s)\times \R\right),
\]
we obtain the Mayer--Vietoris sequence
\[
\cdots \to H_i(\partial W_{\Z})\to H_i(W_\Z)\oplus \bigoplus_{1\leq s \leq r}  H_i(C(V/D_s)\times \R)\to H_i(W_{f,\Z})\to H_{i-1}(\partial W_{\Z})\to \cdots.
\]
We observe that $\partial W_{\Z}$ is homotopy equivalent to a disjoint union of $r$ tori of dimension $(n-2)$ and that the spaces $C(V/D_s)\times \R$ are contractible for $1\leq s\leq r$. Thus, the exactness of the sequence implies that $b_i(W_{\Z})$ is finite if and only if $b_i(W_{f,\Z})$ is finite.
\hfill $\Box$


\section{Middle-dimensional Betti number and finiteness properties}\label{nonfkness}

The aim of this section is to prove Theorem~\ref{theo:nonfiniteness}, showing in particular that when $n=2k$, the group ${\rm ker}(\overline{f}_{\Lambda})$ is not of type $\mathscr{F}_{k}$ for deep enough $\Lambda$, assuming that $\ker (f)$ is of type $\mathscr{F}_{k-1}$. Before proving Theorem~\ref{theo:nonfiniteness} in Section \ref{subsec:proof-nonfiniteness}, we provide a brief discussion of the relevance of $\ell^2$-Betti numbers for the study of finiteness properties in Section \ref{subsec:ell2-betti}. Their use allows for instance to determine the exact finiteness properties of some of the groups considered in~\cite{itmami}. We emphasize however that all the theorems stated in the introduction are proved without using $\ell^{2}$-Betti numbers.

\subsection{The use of $\ell^{2}$-Betti numbers}\label{subsec:ell2-betti}
Informally speaking, the $\ell^2$-Betti numbers of a group $G$ measure the dimension of the (reduced) $G$-invariant $\ell^{2}$-homology of the universal covering space of a $K(G,1)$. When $G$ has a $K(G,1)$ which is a finite complex, the relevant notion of dimension is the {\it von Neumann dimension}. In general, the definition of $\ell^{2}$-Betti numbers is more involved and several approaches are available. See~\cite{chegro1986,gaboriau,kamm-book,luck} for detailed (equivalent) definitions. We will only quickly recall the definition in the proof of Proposition~\ref{observationl2betti} below. The $i$-th $\ell^{2}$-Betti number of $G$ is denoted by $b_{i}^{(2)}(G)$.

The {\it von Neumann algebra} $\mathcal{R}(G)$ of a discrete group $G$ is the algebra of bounded operators on the Hilbert space 
$$\ell^{2}(G)$$
which commute with the left-regular representation of $G$. We will use below that if $M$ is an arbitrary $\mathcal{R}(G)$-module, there is a notion of dimension
$${\rm dim}_{\mathcal{R}(G)} M\in [0,\infty]$$
which was introduced by L\"uck~\cite{luck1998-I}; see also~\cite[Ch. 6]{luck}. When $M$ is finitely generated and projective, one recovers the more classical notion of von Neumann dimension, whose definition can be found e.g. in~\cite[Ch. 1]{luck}. 

The next proposition is an elementary application of a theorem of L\"uck~\cite{luck1998-II}. It already appears in the work of Fisher~\cite{fisher}, in a more general context. We provide the proof for the reader's convenience and we shall comment further on Fisher's work below. In combination with Theorem \ref{thm:Finprops-M8} below, Proposition~\ref{observationl2betti} allows us to determine the precise finiteness properties of the original infinite cyclic coverings considered in \cite{itmami}. 

\begin{prop}\label{observationl2betti} Let $i\ge 1$ be an integer. Let $\Delta$ be a discrete group such that $b_{i}^{(2)}(\Delta)\neq 0$ and let $\phi \colon \Delta \to \Z^{k\geq 1}$ be a surjective homomorphism. Then the kernel of $\phi$ is not of type ${\rm FP}_i(\mathbb{Q})$. 
\end{prop}

\noindent {\it Proof.} Since we are not assuming that the groups under considerations have finite classifying spaces, we must first recall the general definition of $\ell^{2}$-Betti numbers. 

Let $G$ be a group and let $X$ be a $CW$-complex which is a $K(G,1)$. Denote by $\widetilde{X}$ the universal cover of $X$. Let $C_{\ast}(\widetilde{X})$ be the cellular chain complex of $\widetilde{X}$ with $\Z$ coefficients. One now considers the tensor product
$$\mathcal{R}(G)\otimes_{\Z G}C_{\ast}(\widetilde{X}).$$
This gives a chain complex where the underlying modules, and the corresponding homology groups $H_{m}(\mathcal{R}(G)\otimes_{\Z G}C_{\ast}(\widetilde{X}))$ are $\mathcal{R}(G)$-modules. The $m$-th $\ell^2$-Betti number of $G$ is defined as the dimension
$${\rm dim}_{\mathcal{R}(G)}\, H_{m}(\mathcal{R}(G)\otimes_{\Z G}C_{\ast}(\widetilde{X}))$$
of the homology group $H_{m}(\mathcal{R}(G)\otimes_{\Z G}C_{\ast}(\widetilde{X}))$~\cite[\S 4]{luck1998-I}.

Let us now consider the cellular chain complex $C_{\ast}(\widetilde{X},\mathbb{Q})$ of $\widetilde{X}$ with rational coefficients. One sees readily that the complexes
$$\mathcal{R}(G)\otimes_{\Z G}C_{\ast}(\widetilde{X})$$
and $$\mathcal{R}(G)\otimes_{\Q G}C_{\ast}(\widetilde{X},\Q)$$
are isomorphic as $\mathcal{R}(G)$-modules. Hence the $m$-th $\ell^2$-Betti number of $G$ coincides with 
$${\rm dim}_{\mathcal{R}(G)} H_{m}(\mathcal{R}(G)\otimes_{\Q G}C_{\ast}(\widetilde{X},\Q)).$$
The complex $C_{\ast}(\widetilde{X},\Q)$ with the natural augmentation map $C_{0}(\widetilde{X},\Q)\to \Q$ provides a natural free resolution of $\Q$ as a trivial $\Q G$-module. Hence if 
$$\cdots \to F_1 \to F_0 \to \Q\to 0$$
is a free resolution of $\Q$ by $\Q G$-modules, there is a natural homotopy equivalence 
$$C_{\ast}(\widetilde{X},\Q)\to F_{\ast}.$$
After applying $\mathcal{R}(G)\otimes_{\Q G} \cdot$, one sees that the group $H_{m}(\mathcal{R}(G)\otimes_{\Q G} C_{\ast}(\widetilde{X},\Q))$ is isomorphic to $H_{m}(\mathcal{R}(G)\otimes_{\Q G} F_{\ast})$. If $G$ is of type ${\rm FP}_{i}(\Q)$, one can pick a resolution $F_{\ast}$ where the $F_{m}$'s are finitely generated and free for $m\le i$. This implies
$${\rm dim}_{\mathcal{R}(G)} H_{m}(\mathcal{R}(G)\otimes_{\Q G} F_{\ast})<\infty.$$
Hence the $\ell^{2}$-Betti numbers of $G$ are finite in degree $\le i$.

The statement of the proposition now follows by combining this observation with \cite[Th. 3.3 (4)]{luck1998-II}. Indeed, the latter says that if $\Delta$ is a discrete group and if $\psi\colon \Delta\to \Z^{k}$ is a surjective homomorphism, then the finiteness of all the numbers $(b_m^{(2)}(\ker(\psi)))_{0\le m \le i}$ implies the vanishing of all the numbers $(b_m^{(2)}(\Delta))_{0\le m\le i}=0$. See also~\cite[Th\'eor\`eme 6.6]{gaboriau} for a generalization of L\"uck's theorem which also applies here.\hfill $\Box$

As a consequence we obtain:

\begin{cor}\label{cor:observationl2betti}
Let $\Gamma$ be a lattice (uniform or non-uniform) in ${\rm PO}(2n,1)$ or ${\rm PU}(n,1)$ for $n\geq 1$ and let $\phi\colon\Gamma \to \Z^{k\geq 1}$ be a surjective homomorphism. Then $\ker(\phi)$ is not of type ${\rm FP}_{n}(\Q)$.
\end{cor}

Note that the finiteness properties of kernels of homomorphisms from cocompact complex hyperbolic lattices to abelian groups were also studied in~\cite{delzant}, in relation with the BNS invariant. 

\noindent {\emph{Proof of Corollary~\ref{cor:observationl2betti}}.}
By \cite[Theorem 3.3]{kamm-l2-lattices} (see also \cite{chegro} and \cite[Th\'eor\`eme 6.3]{gaboriau}) $b_n^{(2)}(\Gamma)\neq 0$ for every lattice $\Gamma$ in ${\rm PO}(2n,1)$ or ${\rm PU}(n,1)$. Thus, the assertion is an immediate consequence of Proposition \ref{observationl2betti}.
\hfill $\Box$

We provide two applications of Corollary \ref{cor:observationl2betti}.

\begin{exam}\label{ex:itmami}
In \cite{itmami} Italiano, Martelli and Migliorini construct a $2n$-dimensional hyperbolic manifold  $M^{2n}$ together with a homomorphism $f_n\colon \pi_1(M^{2n})\to \Z$ induced by a circle-valued height function for $n\in \left\{2,3,4\right\}$ from duals of Euclidean Gosset polytopes. As a consequence of Corollary \ref{cor:observationl2betti} we obtain that $\ker(f_n)$ is a subgroup of a relatively hyperbolic group that is not of type ${\rm FP}_{n}(\Q)$ (in particular not $\mathscr{F}_{n}$). We emphasize that the fact that the group $\pi_1(M^8)$ is not of type $\mathscr{F}_4$ cannot be obtained from our subsequent arguments for the perturbed map $v\colon M^8\to S^1$ from the introduction, because some of the restrictions of the original map $u$ to the cusps of $M^8$ are null-homotopic, implying that the cusps have infinitely many preimages in the infinite cyclic covering $M^8_{\Z}$.
\end{exam}

We now discuss the following very recent theorem of Fisher~\cite{fisher} that we alluded to in the introduction. To state it we use the notion of {\it Residually Finite Rationally Solvable (short: RFRS)} groups, introduced by Agol (see~\cite{agol} for their definition).

\begin{main} (Fisher)
Let $G$ be a virtually RFRS group of type ${\rm FP}_{n}(\mathbb{Q})$. Then there is a finite index subgroup $H<G$ admitting a homomorphism onto $\mathbb{Z}$ with kernel of type ${\rm FP}_{n}(\mathbb{Q})$ if and only if $b_{i}^{(2)}(G)=0$ for $i=0, \ldots , n$. 
\end{main}

\begin{rem} The proof of the ``only if" part of Fisher's theorem above is given by Proposition~\ref{observationl2betti}. It does not use the RFRS hypothesis. However, Fisher's approach is more general and applies not only to ``usual" $\ell^{2}$-Betti numbers but also to a class of generalized Betti numbers with skew-field coefficients.
\end{rem}

Consider now a torsionfree cocompact lattice $\Gamma < {\rm PO}(2k,1)$ which is cubulable. One can take $\Gamma$ to be any torsionfree arithmetic lattice of the {\it simplest type}, thanks to the work of Bergeron, Haglund and Wise~\cite{bhw}. Being hyperbolic and cubulable, such a lattice is virtually special due to Agol's Theorem~\cite[Theorem 1.1]{Agol13}. In particular $\Gamma$ must be virtually RFRS. The $\ell^{2}$-Betti numbers of any finite index subgroup $\Gamma'$ of $\Gamma$ vanish up to degree $k-1$ and $b_{k}^{(2)}(\Gamma')\neq 0$ \cite[Theorem 3.3]{kamm-l2-lattices}. Applying Fisher's theorem we can thus find a finite index subgroup $\Gamma'<\Gamma$ and a surjective homomorphism $f \colon \Gamma'\to \mathbb{Z}$ such that $\ker(f)$ is of type ${\rm FP}_{k-1}(\mathbb{Q})$. The non-vanishing of $b_{k}^{(2)}(\Gamma')$ together with Proposition~\ref{observationl2betti} implies that $\ker(f)$ is not of type ${\rm FP}_{k}(\Q)$. This proves:

\begin{prop} Let $\Gamma < {\rm PO}(2k,1)$ be a torsionfree cocompact cubulable lattice. Then there exists a finite index subgroup $\Gamma'< \Gamma$ and a surjective homomorphism $f \colon \Gamma'\to \mathbb{Z}$ such that $\ker(f)$ is of type ${\rm FP}_{k-1}(\mathbb{Q})$ but not of type ${\rm FP}_{k}(\mathbb{Q})$.
\end{prop}

For every integer $k\ge 1$ this produces plenty of subgroups of hyperbolic groups which are ${\rm FP}_{k-1}(\mathbb{Q})$ but not ${\rm FP}_{k}(\mathbb{Q})$.

In the next section, we will prove Theorem~\ref{theo:nonfiniteness} using arguments involving the homology with rational coefficients of the cyclic covering spaces associated to the characters under consideration. An alternative approach to prove Theorem~\ref{theo:nonfiniteness} would be to compute the $\ell^2$-Betti numbers of $\Lambda/N_{f}(\Lambda)$ 
and to appeal to L\"uck's theorem~\cite[Th. 3.3]{luck1998-II}. However, as far as we are aware, the $\ell^{2}$-Betti numbers of the group $\Lambda/N_{f}(\Lambda)$ have not been computed. More generally, it would be interesting to study the behavior of $\ell^{2}$-Betti numbers under Dehn fillings. We spell out one concrete question in this direction. 

\begin{question} Let $\Gamma < {\rm PO}(2k,1)$ be a torsionfree nonuniform lattice with purely unipotent parabolic subgroups. Let $\Lambda < \Gamma$ be a finite index subgroup. Let $(H_{s})_{1\le s \le r}$ be pairwise nonconjugate maximal parabolic subgroups of $\Lambda$, representing all conjugacy classes of maximal parabolic subgroups. For every $s\in \left\{1, \ldots , r\right\}$, pick a subgroup $A_s < H_s$. 

Is it true that the $k$-th $\ell^{2}$-Betti number of the quotient group
$$\Lambda/\langle \langle \bigcup_{1\le i \le N} A_{i}\rangle \rangle$$
of $\Lambda$ by the normal closure of $\bigcup_{1\le i \le N} A_{i}$ is nonzero if $\Lambda$ is deep enough? 
\end{question}

This seems to be related to L\"uck's Approximation Conjecture~\cite[Ch. 13]{luck}. We would like to thank Jean Raimbault for pointing out this connection to us.

\subsection{Middle-dimensional Betti number and infinite cyclic coverings}\label{subsec:proof-nonfiniteness}

We resume with the notation of Section \ref{dehnfillingfiniteness}. In addition we now assume that $n=2k$ is even. By Corollary~\ref{cor:observationl2betti}, $\ker (f)$ is not ${\rm FP}_{k}(\Q)$. We will prove the following result, which does not require the use of $\ell^{2}$-Betti numbers, but where instead we assume that $b_i(\ker(f))$ is finite for $i\le k-1$. This is for instance the case when $\ker (f)$ is of type $\mathscr{F}_{k-1}$.

\begin{prop}\label{prop:finitude} Assume that $b_{i}(\ker(f))$ is finite for $i\le k-1$. Then the group $H_{k}(\ker (f))$ is infinite dimensional. In particular, the group $\ker(f)$ is not of type $\mathscr{F}_{k}$.
\end{prop} 

This proposition implies the first statement in Theorem~\ref{theo:nonfiniteness}. Note that since the Euler characteristic of $\mathbb{H}^{2k}/\Gamma$ is nonzero~\cite{kz}, it follows from Milnor \cite[Assertion 6]{milnor} that $\mathbb{H}^{2k}/\ker(f)$ has {\it some} Betti number which is infinite. The proof of Proposition~\ref{prop:finitude} consists in proving that the $k$-th Betti number of this space is infinite (and only this one).

We will require the following consequence of Poincar\'e-Lefschetz duality:
\begin{lemma}\label{lem:Poincare-Lefschetz}
 Let $(M,\partial M)$ be an $n$-dimensional compact oriented manifold with boundary and let $i\in \{1, \ldots , n\}$. Then 
 \[
   b_{n-i}(M)\leq b_{i-1}(\partial M) + b_i(M).
 \]
\end{lemma}
\noindent {\it Proof.}
By Poincar\'e Lefschetz duality and the Universal Coefficient Theorem there are isomorphisms
\[
 H_i(M)\cong H^{n-i}(M,\partial M)\cong H_{n-i}(M,\partial M)
\]
for $0\leq i \leq n$.
From the long exact sequence
\[
\cdots \rightarrow H_{n-i}(\partial M)\rightarrow H_{n-i}(M)\rightarrow H_{n-i}(M,\partial M)\rightarrow H_{n-i-1}(\partial M)\rightarrow \cdots
\]
of the pair $(M,\del M)$ in homology, we deduce that
\[
 b_{n-i}(M)\leq b_{n-i}(\partial M)+b_{n-i}(M,\partial M)\leq b_{n-i}(\partial M)+b_i(M) = b_{i-1}(\partial M)+ b_{i}(M),
\]
where for the last equality we use Poincar\'e duality for the closed manifold $\partial M$.
\hfill $\Box$

\noindent {\it Proof of Proposition \ref{prop:finitude} .} We write $M_{\Z}=\mathbb{H}^{2k}/\ker(f)$ and $M_{\ell}=\mathbb{H}^{2k}/N_{\ell}$, where $N_{\ell}$ is the kernel of the morphism 
$$\Gamma \to \Z/\ell\Z$$
obtained by reducing $f$ mod $\ell$. Let $T \colon M_{\Z}\to M_{\Z}$ be a generator of the group of deck transformations of the cyclic covering space $M_{\Z}\to \mathbb{H}^{2k}/\Gamma$. Let $t \colon H_{\ast}(M_{\Z})\to H_{\ast}(M_{\Z})$ be the transformation induced by $T$ on homology. By applying Milnor's long exact sequence~\cite[p. 118]{milnor} to the cyclic covering $M_{\Z}\to M_{\ell}$, we obtain the following short exact sequence:
\begin{equation}\label{exactsequence}
0\to H_{i}(M_{\Z})/(t^{\ell}-1)H_{i}(M_{\Z})\to H_{i}(M_{\ell})\to \ker(t^{\ell}-1 \colon H_{i-1}(M_{\Z})\to H_{i-1}(M_{\Z}))\to 0.
\end{equation} 
For $i\le k-1$ the vector spaces appearing on the left and right in the above short exact sequence have finite dimension bounded by $b_{i}(M_{\Z})$ and $b_{i-1}(M_{\Z})$ respectively. We thus obtain
\begin{equation}
\label{eqn:Betti-M-ell}
b_{i}(M_{\ell})\le b_{i}(M_{\Z})+b_{i-1}(M_{\Z})
\end{equation}
for $i\le k-1$ and all $\ell\ge 1$. 

Since the restriction of $f$ to each of the finitely many cusps of $M$ is non-trivial, the number of boundary tori in $M_\ell$ is independent of $\ell$. Therefore, for all $i$, the Betti numbers $b_{i-1}(\partial M_{\ell})$ are also independent of $\ell$.
By Lemma \ref{lem:Poincare-Lefschetz} we have that $b_{2k-i}(M_{\ell})\leq b_{i-1}(\partial M_{\ell})+b_i(M_{\ell})$ for $i\leq k-1$. Thus, we deduce from \eqref{eqn:Betti-M-ell} that all Betti numbers $b_{i}(M_{\ell})$, for $i\in \{0, \ldots, 2k\}$ with $i\neq k$, are uniformly bounded above (independently of $\ell$). 

Since the Euler characteristic of $M_{\ell}$ is equal to $\ell$ times that of $M$ (which is nonzero~\cite{kz}), we obtain that $b_{k}(M_{\ell})$ grows roughly linearly with $\ell$. In particular
$$b_{k}(M_{\ell})\to \infty$$
as $\ell$ goes to $\infty$. Combining Equation~\eqref{exactsequence}, this time for $i=k$ with the fact that $b_{k-1}(M_{\Z})$ is finite, we deduce that the codimension of the image of the natural map 
$$H_{k}(M_{\Z})\to H_{k}(M_{\ell})$$
is bounded above uniformly in $\ell$. This implies that $H_{k}(M_{\Z})$ is infinite dimensional and completes the proof.\hfill $\Box$

As a consequence we obtain the following corollary. It corresponds to the second assertion in Theorem~\ref{theo:nonfiniteness}, thus completing its proof.

\begin{cor} We keep the asumptions from Proposition~\ref{prop:finitude}. Let $\Lambda<\Gamma$ be a deep enough finite index subgroup. Then the $k$-th Betti number of the group $\ker (\overline{f}_{\Lambda})$ is infinite. In particular $\ker(\overline{f}_{\Lambda})$ is not of type $\mathscr{F}_{k}$. 
\end{cor}

\noindent {\it Proof.} 
This is an immediate consequence of Proposition \ref{prop:finitude} and Proposition \ref{prop:Betti-finiteness}. 
\hfill $\Box$


\section{The manifold $M^8$}\label{sec:m8}

In this section we describe the hyperbolic $8$-manifold $M^8$ first studied in~\cite{itmami} and prove Theorem~\ref{goodmaptothecircle}. We will use the letter $M$ to denote a general hyperbolic manifold and will write $M^{8}$ when refering to the specific example built in~\cite{itmami}. 

\subsection{Colourings}

Let $P\subset \Xs^n$ be a finite-volume right-angled polytope, with $\Xs^n = \R^n$ or $\Hy^n$. A \emph{facet} is a codimension-1 face of $P$. A $c$-\emph{colouring} on $P$ is the assignment of a colour $\lambda(F) \in \{1,\ldots, c\}$ to each facet $F$ of $P$, such that adjacent facets have distinct colours. We always suppose that each colour in the palette $\{1,\ldots, c\}$ is assigned to at least one facet.

It is a standard fact, which goes back to \cite{V}, that a colouring produces a manifold $M$ with the same geometry as $\Xs^n$, tessellated into $2^c$ copies of $P$. The construction goes as follows.
The reflection group $\Gamma_0$ of $P$ is generated by the reflections $r_F$ along the facets $F$ of $P$. A colouring defines a homomorphism $\Gamma_0 \to \Z_2^c$ which sends $r_F$ to $e_{\lambda(F)}$, the generator of the $\lambda(F)$-th factor of $\Z_2^c$. The kernel $\Gamma$ of this homomorphism is torsion-free, and hence defines a (flat or hyperbolic) manifold $M = \Xs^n/\Gamma$.  We get an orbifold-covering $M\to P$ and $M$ is tessellated into $2^c$ copies of $P$.

Here is a concrete description of the tessellation. For every $v\in \Z_2^c$ we pick a copy $P_v$ of $P$. We glue these $2^c$ copies $\{P_v\}$ as follows: we identify every facet $F$ of $P_v$ with the same facet $F$ (via the identity map) of the polytope $P_{v+e_i}$ with $i = \lambda(F)$. The result is the tessellation of $M$.

If $P$ is compact, then $M$ also is. If $P\subset \Hy^n$ has some ideal vertices, the manifold $M$ is finite-volume and cusped. Let $v$ be an ideal vertex of $P$. A link $\lk(v)$ of $v$ in $P$ is a right-angled Euclidean parallelepiped, which inherits a colouring (still denoted by $\lambda$) in the obvious way: every facet $F$ of $\lk(v)$ is contained in a unique facet $F'$ of $P$ and we set $\lambda(F) = \lambda(F')$. The colouring $\lambda$ on $\lk(v)$ produces an abstract flat compact $(n-1)$-manifold $T$ which is topologically a $(n-1)$-torus \cite[Proposition 7]{itmami} tessellated into $2^{c'}$ copies of $\lk(v)$, where $c'\leq c$ is the number of distinct colours inherited by $\lk(v)$.

The preimage of $\lk(v)$ along the orbifold covering $M \to P$ consists of $2^{c-c'}$ copies of $T$. In particular all the cusps of $M$ are toric. Some instructing examples are exposed in \cite{itmami}.

\subsection{The manifold $M^8$}
The 8-dimensional cusped hyperbolic manifold $M^8$ is constructed in \cite{itmami} by taking the right-angled polytope $P=P^8 \subset \Hy^8$ dual to the Euclidean Gosset polytope $4_{21}$, whose vertices are the 240 non-trivial elements of $E_8$ of smallest norm. The polytope $P^8$ has 240 facets, 2160 ideal vertices, and 17280 finite vertices. Its isometry group acts transitively on the facets.

Using octonions, the authors defined in \cite{itmami} an extremely symmetric 15-colouring for $P^8$, where each colour is assigned to 16 distinct pairwise non-incident facets. It is also shown in \cite[Proposition 10]{itmami} that this colouring is minimal, that is, $P^8$ cannot be 14-coloured.

As shown in \cite{itmami} there are two types of ideal vertices $v$ with respect to this colouring: in total there are 1920 ideal vertices of the first type and 240 of the second type. The link $\lk(v)$ of a vertex $v$ of the first type is a 14-coloured 7-cube, while the link of a vertex of the second type is a 7-coloured 7-cube. Each vertex $v$ of the first type produces $2^{15-7} = 256$ cusps and each vertex of the second type produces $2^{15-14} = 2$ cusps. So the manifold $M^8$ has a total of $240\cdot 256 + 1920 \cdot 2 = 65280$ cusps.

\subsection{The combinatorial game}\label{subsec:JNW-game}
The combinatorial game of Jankiewicz, Norin and Wise \cite{JNW} is then applied in \cite{itmami} to construct a nice map $f\colon M^8 \to S^1$. We briefly explain the rules of the game, since we will need them below. We introduce here a slightly modified version that is more adapted to our purposes and has already appeared in \cite{itmami2}.

We start with a right-angled polytope $P$, equipped with a colouring $\lambda$ that produces a manifold $M$. A \emph{state} for $P$ is an assignment of a status I or O to each facet. The letters stand for In and Out. \emph{A set of moves} is a partition of the colour palette $\{1,\ldots, c\}$. The individual sets of the partition are called \emph{moves}.

The tessellation of $M$ into $2^c$ copies of $P$ is dual to a cubulation $C$, as explained in \cite[Section 2.2]{itmami}.  
The cubulation $C$ has $2^c$ vertices dual to the copies of $P$ in $M$, its edges are dual to their facets, and so on. If $P$ has some ideal vertices, the cubulation is a compact object that embeds naturally in $M$ as a spine: its complement $M\setminus C$ consists of the cusps. In particular $M$ deformation retracts onto $C$.
As explained in \cite[Section 1.6]{itmami2}, the state and the set of moves determine an orientation of all the edges of the dual cubulation $C$. This can be done inductively as follows. The vertices of $C$ are naturally identified with $\Z_2^c$. We start with the vertex $0\in \Z_2^c$: the edges adjacent to $0$ are dual to the facets of $P$, and an edge $s$ connects $0$ and $e_i$ if the dual facet $F$ is coloured by $i$. (There are multiple edges with the same endpoints if there are multiple facets having the same colour.) We orient the edge $s$ towards $0$ (inward) if $F$ has status I, and towards $e_i$ (outward) if it has status O. Now suppose that the edges adjacent to the vertex $v\in \Z_2^c$ are all oriented, and we must decide how to orient the edges adjacent to $v+e_i$ for some $i$. We note that the edges connecting $v$ and $v+e_i$ are dual to the facets of $P$ coloured with $i$.

The edges adjacent to $v+e_i$ are in natural 1-1 correspondence with those adjacent to $v$, since they are both in 1-1 correspondence with the facets of $P$. The rule is the following: when we pass from $v$ to $v+e_i$, we invert the orientation (as seen from $v$ and $v+e_i$ respectively) of all the edges whose dual facet in $P$ has a colour that lies in the same move containing $i$. We can prove easily that these rules yield a well-defined orientation on all the edges of $C$.

\begin{figure}
 \begin{center}
  \includegraphics[width = 8 cm]{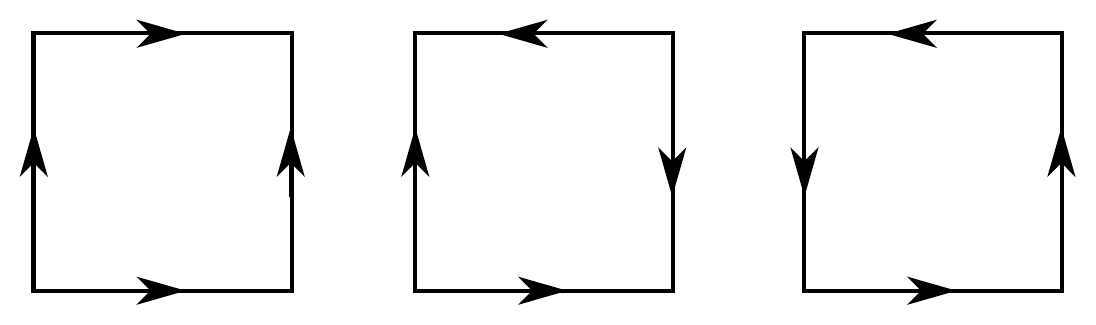}
 \end{center}
 \caption{The possible configurations that may arise on a square of the  cubulation $C$.}
 \label{squares:fig}
\end{figure}

We have oriented all the edges of the cubulation $C$, and we now examine the possible configurations at the squares. Every square of $C$ is dual to a codimension-2 face of $P$ that separates two facets $F_1$ and $F_2$ of $P$ with distinct colours $i_1$ and $i_2$. There are three possibilities:
\begin{enumerate}
    \item
If $i_1$ and $i_2$ belong to two distinct moves, the edges are oriented as in Figure \ref{squares:fig}-(left); in this case we say that the pair $F_1, F_2$ is \emph{very good};
\item If $i_1$ and $i_2$ belong to the same move, and $F_1$ and $F_2$ have the same status, the edges are oriented as in Figure \ref{squares:fig}-(centre); we say that the pair $F_1, F_2$ is \emph{good};
\item if $i_1$ and $i_2$ belong to the same move and $F_1$ and $F_2$ have opposite stati, the edges are oriented as in Figure \ref{squares:fig}-(right) and we say that the pair $F_1, F_2$ is \emph{bad}. 
\end{enumerate}

If every pair of adjacent facets in $P$ is very good, we say that the orientation of the edges is \emph{coherent}, and as explained in \cite{itmami} we can define a \emph{diagonal map} $f\colon C \to \R/\Z = S^1$ which is a Morse function in the sense of Bestvina and Brady \cite{BB}. We can combine this map with the deformation retraction of $M$ onto $C$ to get a map $f\colon M \to S^1$.

If every pair of adjacent facets in $P$ is either good or very good, we do not get a diagonal map as above, but if we assign the same number 1 to each oriented edge we still get a 1-cocycle: indeed we see from Figure \ref{squares:fig} that the contribution of the four edges sums to zero. This is a crucial point employed in \cite{itmami2}.

If the set of moves is just the discrete partition
$$\{1\}, \ldots, \{c\}$$
we deduce that distinct colours always belong to distinct moves, and hence all the squares are very good. This is the set of moves that was used implicitly in all the examples of \cite{itmami}. We will need a different kind of set of moves at some point below.

\subsection{The abelian cover $M^8_{\Z}$}
Using octonions, the authors have defined in \cite{itmami} a particularly symmetric state $s$ for $P^8$. The state is \emph{balanced}, in the sense that for each 16 facets sharing the same colour, half of them receive the status O and the other half receive the status I. They used as a set of moves the discrete partition, and hence they got a Morse function $f\colon C \to S^1$ inducing a map $f\colon M^8 \to S^1$.

The authors have then used a program written in Sage and available from \cite{code} to analyse the ascending and descending links of all the vertices of the cubulation, see \cite{BB} for the terminology. The output of the program shows that all the ascending and descending links are simply connected, and hence by \cite[Theorem 4.1]{BB} the subgroup $\ker (f)$ is finitely presented. Therefore the infinite cyclic manifold cover $M^8_{\Z}\to M^8$ determined by $\ker (f)$ is a complete (infinite volume, geometrically infinite) hyperbolic 8-manifold with finitely presented fundamental group $\ker (f)$. However $\ker (f)$ is not of type $\mathscr{F}_7$ because $b_7(M^8_{\Z})=\infty$, see \cite[Theorem 23]{itmami}.

In fact, more can be proved.

\begin{main}\label{thm:Finprops-M8}
The group $\ker (f)$ is of type $\mathscr{F}_3$ and not of type $\mathscr{F}_4$.
\end{main}
\noindent {\it Proof.} 
The same code from \cite{code} shows that every ascending and descending link is also homologicaly 2-connected. From \cite[Theorem 4.1, Lemma 3.7]{BB} we deduce that $\ker (f)$ is ${\rm FP}_3$, which together with finite presentability implies $\mathscr{F}_3$ \cite[Proof of Theorem 7.1]{brown}. On the other hand Example \ref{ex:itmami} shows that $\ker(f)$ is not of type $\mathscr{F}_4$.
\hfill $\Box$ 

\subsection{Cusp homologies}\label{subsec:cusphom}
We now turn back to the more general setting of a right-angled polytope $P \subset \Hy^n$ equipped with a colouring $\lambda$ that produces a manifold $M$. As we said above, every ideal vertex $v$ of $P$ gives rise to some cusps in $M$ that lie above $v$. That is, the preimage in $M$ of a link $\lk(v)$ of an ideal vertex $v$ in $P$ consists of a union of disjoint flat $(n-1)$-tori. The orbifold covering $M\to P$ is regular, so the deck transformation group acts transitively on these tori. We pick one such cusp section $\iota \colon T \hookrightarrow M$ and we are interested in the induced map $\iota^*\colon H^1(M, \R) \to H^1(T, \R)$ in cohomology.

\begin{prop} \label{C:prop}
Suppose that there are two opposite facets $F_1, F_2$ of $\lk(v)$ such that one of the following holds:
\begin{enumerate}
    \item The facets $F_1$ and $F_2$ have the same colour, or
    \item The facets $F_1$ and $F_2$ have distinct colours $i_1$ and $i_2$, and they lie in distinct components of the subset $X\subset \partial P$ consisting of all the facets coloured with either $i_1$ or $i_2$. 
\end{enumerate}
Then the map $\iota^*\colon H^1(M, \R) \to H^1(T, \R)$ is non-trivial. If this holds for every pair of opposite facets of $\lk(v)$, the map $\iota^*$ is surjective.
\end{prop}
\noindent {\it Proof.}
Let $F_{1,i}, F_{2,i}$ denote the pairs of opposite facets in $\lk(v)$, for $i=1,\ldots,n-1$. Note that only opposite facets may share the same colour. Let $s_i\subset \lk(v)$ be the segment connecting the centers of $F_{1,i}$ and $F_{2,i}$. The preimage of $F_{1,i}\cup F_{2,i}$ in $T$ consists of either two or four parallel facets, depending on whether $F_{1,i}$ and $F_{2,i}$ share the same colour or not. The preimage of $s_i$ in $T$ is a geodesic loop $\gamma_i$ orthogonal to these, consisting of either two or four copies of $s_i$. We get $n-1$ geodesic loops $\gamma_1,\ldots,\gamma_{n-1} \subset T$ that are pairwise orthogonal and generate $H_1(T, \R)$. Let $\alpha_1, \ldots, \alpha_{n-1} \in H^1(T, \R)$ be the dual basis.

We now show that if the pair $F_{1,j}, F_{2,j}$ satisfies either (1) or (2) then $\alpha_j$ lies in the image of $\iota^*$. This will conclude the proof. 

In case (1), fix a state $s$ for $P$ such that $F_{1,j}$ and $F_{2,j}$ have opposite stati, while all the other facets of $\lk(v)$ have status I. Use the discrete partition $\{1\}, \ldots, \{c\}$ as a set of moves. Every pair of facets in $P$ is very good, so we get an orientation of the edges of the dual cubulation $C$ which induces a 1-cocycle and therefore a class $\alpha \in H^1(M, \R)$. The image $\iota^*(\alpha) \in H^1(T, \R)$ sends $\gamma_j$ to $\pm 2$ and $\gamma_i$ to zero for all $i\neq j$. This holds because each $\gamma_i$ is isotopic to the cycle of two or four edges in the dual cubulation dual to the lifts of $F_{1,i}$ and $F_{2,i}$, and two consecutive edges are oriented in the same direction with respect to their common endpoint $v$ (both inward or both outward with respect to $v$) if and only if the two stati of $F_{1,i}$ and $F_{2,i}$ coincide. Therefore $\iota^*(\alpha) = \pm 2 \alpha_j$ and we are done.

In case (2), fix a state $s$ for $P$ such that $F_{1,j}$ and $F_{2,j}$ have opposite stati, and in each connected component of $X$ all the facets have the same status (we can do this because $F_{1,j}$ and $F_{2,j}$ lie in distinct connected components). We also suppose as above that all the other facets of $\lk(v)$ have status I. We suppose that the colours of $F_{1,j}$ and $F_{2,j}$ are 1 and 2 for simplicity of notation and as a set of moves we pick $\{1,2\}, \{3\}, \{4\}, \ldots, \{c\}$. By construction every pair of adjacent facets in $P$ is either very good or good. Therefore the resulting orientation of the edges of the cubulation yields a cocycle $\alpha$. As in (1) we conclude that $\iota^*(\alpha)$ sends $\gamma_j$ to $\pm 4$ and $\gamma_i$ to 0 for all $i\neq j$. Therefore $\iota^*(\alpha) = \pm 4\alpha_j$.
\hfill $\Box$ 

We remark that it is also possible to prove the proposition using the work of Choi and Park on the cohomology of $M$, see \cite{CP}. We also note that there are some cases where $\iota^*$ is indeed trivial: pick a $2n$-gon $P\subset \Hy^2$ with $2n-1$ right-angled real vertices and one ideal vertex, and colour $P$ with 2 colours. None of the hypotheses of the proposition are fulfilled, and indeed, in this case, $M$ is an orientable surface with a single cusp, and hence $\iota^*$ is trivial  (here $T=S^1$).

We will need here the following consequence.

\begin{prop}\label{C:prop-M8}
The cohomology group $H^1(M^8,\R)$ has dimension 365. For every cusp section $T$ the natural map $H^1(M^8, \R) \to H^1(T, \R)$ is surjective.
\end{prop}
\noindent {\it Proof.} 
The number 365 was calculated in \cite{itmami} using a formula of Choi -- Park \cite{CP}.
The 7-torus cusp section $T$ projects to a 7-cube $\lk(v)$ of some ideal vertex $v$ of $P^8$. There are two types of 7-cubes: for the 7-coloured ones, two opposite facets share the same colours and then Proposition \ref{C:prop} applies; for the 14-coloured ones, two opposite facets have distinct colours, but we have checked using Sage that the property stated in Proposition \ref{C:prop}-(2) holds for all pairs of opposite facets. The code is available from \cite{code}.
\hfill $\Box$

\subsection{Constructing perturbations}

We start this section with a short review of the theory of BNSR-invariants that will be useful to perturb certain characters while keeping control on the finiteness properties of their kernels. Let $G$ be an arbitrary finitely generated group. A $\emph{character}$ of $G$ is a homomorphism $\chi \colon G \to \R$. This is the same as an element $\chi\in H^1(G,\R)$. We introduce the {\it character sphere} of equivalence classes of characters
\[
S(G) = (H^1(G,\R)\setminus \left\{0\right\})/\sim,
\]
where we call two characters $\chi_1$ and $\chi_2$ \emph{equivalent} and write $\chi_1\sim \chi_2$ if there is $\lambda \in \R_{>0}$ with $\lambda \cdot \chi_1 = \chi_2$. We denote by $[\chi]$ the equivalence class of a character $\chi$ in $S(G)$. We call $[\chi]$ {\emph{rational}} if it has a representative $\chi\colon G \to \Q$. Note that the rational equivalence classes of characters form a dense subset of $S(G)$.

We equip $S(G)$ with the topology induced by a choice of auxiliary norm on the real vector space $H^1(G,\R)$. Then $S(G)$ is a sphere of dimension $b_1(G)-1$.

In \cite{Renz-thesis}, Renz proves that the finiteness properties of the kernels of characters are determined by a descending sequence of subsets ${}^*\Sigma ^i(G)$ of $S(G)$
\[
S(G)= {}^* \Sigma^0(G) \supseteq {}^*\Sigma^1(G) \supseteq {}^*\Sigma^2(G) \supseteq \dots,
\]
called the geometric (BNSR-)invariants of $G$. We refer to \cite{Renz-thesis} for their precise definition.

More precisely Renz proves \cite[Corollary AC]{Renz-thesis} (see also \cite[p.477, p.481]{Renz-sing}, \cite[Remark 6.5]{BieriRenz} and~\cite{bieristrebel}):
\begin{main}\label{thm:renz}
Let $G$ be a group of type $\mathscr{F}_k$ and let $S(G)$ be its character sphere. Then there is a descending chain of open subsets $\left\{{}^*\Sigma^i(G)\right\}_{i\geq 0}$ of $S(G)$, such that the kernel of a rational character $\chi \colon G\to \Q$ is of type $\mathscr{F}_i$ for some $0\leq i \leq k$ if and only if $[\chi], [-\chi]\in {}^*\Sigma^i(G)$. In particular, for $0\leq i \leq k$, the condition that the kernel $\ker(\chi)$ of a rational equivalence class of characters $[\chi]$ is of type $\mathscr{F}_i$ is an open condition in $S(G)$.
\end{main}

We now return to the setting of Section \ref{subsec:cusphom} and consider a hyperbolic manifold $M$ with fundamental group $\Gamma$. 
The main result of this section is the following. 

\begin{prop}\label{prop:perturb}
Let $M$ and $\Gamma$ be as above and let $\chi\colon\Gamma \to \Z$ be a character with kernel of type $\mathscr{F}_k$ for some $k\geq 0$. Let $T_1,\ldots, T_r$ be cusp sections of the $r$ cusps of $M$. Assume that for all $1\leq s \leq r$ the morphism $\iota_s^{\ast}\colon H^1(M,\R)\to H^1(T_s,\R)$ induced by the inclusion $\iota_s\colon T_s\hookrightarrow M$ is non-trivial. 

Then the set of all rational characters $U\subset H^1(M,\Q)$ with the properties that
\begin{enumerate}
    \item $\forall \mu \in U$ the kernel $ker(\mu)$ is of type $\mathscr{F}_k$, and
    \item $\forall \mu \in U$ the restriction $\mu|_{\pi_1(T_s)}$ to every cusp $T_s$, $1\leq s \leq r$, is non-trivial,
\end{enumerate}
is a non-empty open subset of $H^1(M,\Q)\setminus \left\{0\right\}$ with respect to the subspace topology induced by $H^1(M,\R)$. More precisely, $U$ can be obtained by intersecting $H^1(M,\Q)\setminus \left\{0\right\}$ with an infinite cone over a non-empty open subset of $S(\Gamma)$.
\end{prop}
\noindent {\emph{Proof.}}
 First observe that the set
 \[
  U_0:= H^1(M,\Q)\setminus \left( \bigcup_{1\leq s\leq r} (\iota_s^{\ast})^{-1}(0)\right)
 \]
 of all rational characters $\mu\colon\Gamma \to \Q$ which satisfy condition (2) is a dense open subset of $H^1(M,\Q)$. 
 
Since $\Gamma$ is of type $\mathscr{F}_k$, the existence of a character $\chi$ with kernel of type $\mathscr{F}_k$ and Theorem \ref{thm:renz} imply that the subset of $S(\Gamma)$ consisting of equivalence classes of characters with kernel of type $\mathscr{F}_k$ is non-empty and open. Combining this with the above description of $U_0$ yields the desired conclusion.
\hfill $\Box$

Proposition \ref{prop:perturb} allows us to perturb characters $\Gamma\to \Z$ so that they induce non-trivial homomorphisms on all parabolic subgroups of $\Gamma$. Thus, we can now complete the

\noindent {\emph{Proof of Theorem \ref{goodmaptothecircle}.}}
 By Theorem \ref{thm:Finprops-M8} and Proposition \ref{C:prop-M8} the manifold $M^8$ together with the character $u_*\colon \Gamma \to \Z$ satisfies the assumptions of Proposition \ref{prop:perturb} for $k=3$. Thus there is a character $\chi\colon \Gamma \to \Z$ satisfying properties (1) and (2) of Proposition \ref{prop:perturb}. Since every character on $\Gamma$ is induced by a continuous map $v\colon M \to S^1$, this concludes the proof.
\hfill $\Box$

\section{An infinite family of examples}\label{sec:infinitefamily}
In this section we will explain how to obtain an infinite family of pairwise non-isomorphic hyperbolic groups with subgroups of type $\mathscr{F}_3$ and not $\mathscr{F}_4$ from our construction. 

Let $M= \Hy^n/\Gamma$ be a finite-volume hyperbolic $n$-manifold with $r$ toric cusps, with sections $T_1,\ldots, T_r$.

\begin{prop}\label{prop:char-systoles}
 Assume that for every cusp section $T_s$ the natural map $H^1(M,\R)\to H^1(T_s,\R)$ is surjective and that there is a character $\chi\colon\Gamma \to \Z$ with kernel of type $\mathscr{F}_k$ for some $k\geq 0$. Then there is a sequence of characters $\mu_n\colon\Gamma \to \Z$, $n\in \N$, with the following properties: 
 \begin{enumerate}
     \item for every $n\in \N$ the kernel $\ker(\mu_n)$ is of type $\mathscr{F}_k$,
     \item for every $n\in \N$ and every cusp section $T_s$ of $M$ the systole of the natural subtorus associated to the subgroup $\ker(\mu_n)\cap \pi_{1}(T_{s}) < \pi_{1}(T_s)$ equipped with the flat metric induced by $T_s$ is $\geq n$.
 \end{enumerate}
\end{prop}
\noindent {\emph{Proof.}}
 Let $V\subset H^1(M,\Q)\setminus \left\{0\right\}$ be the open cone consisting of rational characters whose kernel is of type $\mathscr{F}_k$.  
 
 Note that for every $s\in \left\{1,\dots,r\right\}$ and for every $n\in \N$ the set $Z_{n,s}\subset \pi_1(T_s)$ of non-trivial homotopy classes of loops in $T_s$ that are homotopic to a loop of length $\leq n$ is finite. The surjectivity of $\iota_s^*\colon H^1(M,\R)\to H^1(T_s,\R)$ thus allows us to choose characters $\mu_{n,s}\colon\Gamma \to \Z$ such that $Z_{n,s}\cap \ker(\mu_{n,s})=\emptyset$ for $1\leq s \leq r$. Using the openness of $V$ and a simple induction on $s$ we can choose arbitrarily small rational numbers $\lambda_{n,s}\in \Q\setminus \left\{0\right\}$ such that the character 
 \[
 \mu_n := \chi + \sum_{s=1}^r \lambda_{n,s} \mu_{n,s}\colon \Gamma \to \Q
 \]
 is in $V$ and satisfies $\ker(\mu_n)\cap Z_{n,s}=\emptyset$ for $1\leq s \leq r$. Replacing $\mu_n$ by a representative of its equivalence class with values in $\Z$ completes the proof.
\hfill $\Box$

As a consequence we obtain:
\begin{prop}\label{prop:FinProps-inf-family}
With the same assumptions as in Proposition \ref{prop:char-systoles}, assume now that ${\rm dim}(M) =2k$. Then there is an infinite family of pairwise non-isomorphic hyperbolic groups, each of which has a subgroup which is of type $\mathscr{F}_{k-1}$ and not of type $\mathscr{F}_{k}$.
\end{prop}

\noindent {\emph{Proof.}}
 Let $\mu_n\colon\Gamma \to \Z$ be a sequence of characters as in Proposition \ref{prop:char-systoles}. By construction the intersection 
 \begin{equation}
 \label{eqn:triv-intersection}
 \bigcap_{n\in \N} \ker \left(\mu_n\circ \iota_{s,\ast}\colon \pi_1(T_s)\to \Z\right) =\left\{1\right\}
 \end{equation}
 is trivial for $1\leq s \leq r$. As described in Section \ref{dehnfi} we can use the methods of Fujiwara and Manning \cite{fujiwaramanning} to construct negatively curved Dehn fillings $\overline{M}_n$ for all characters $\mu_n$ with $n\geq 2\pi$ by coning off a family of translates over $S^1$ of the subtori of the $T_s$ associated to the kernels $\ker(\mu_n)$. Note that for $n\geq 2\pi$ all fillings are $2\pi$-fillings and therefore a passage to a finite index subgroup of $\Gamma$ is not required before performing the Dehn filling in order to apply Theorem \ref{theo:fujiwaramanning} and obtain that $\overline{M}_n$ is ${\rm CAT}(-1)$.
 
 It then follows from \eqref{eqn:triv-intersection} and \cite[Proposition 2.12 and Theorem 2.13]{fujiwaramanning} that the family $\left\{\pi_1(\overline{M}_n)\mid n\in \N \right\}$ contains an infinite family of pairwise non-isomorphic hyperbolic groups. By combining Theorem \ref{fkispreserved} and Theorem \ref{theo:nonfiniteness}, we obtain that the induced homomorphisms $\overline{\mu_n}\colon \pi_1(\overline{M}_n)\to \Z$ all have kernel of type $\mathcal{F}_{k-1}$ and not of type $\mathcal{F}_{k}$, thus completing the proof.
\hfill $\Box$

We can now conclude:

\noindent {\emph{Proof of Theorem \ref{ssgrofhyperbolicgroup}}}.
By Theorem \ref{thm:Finprops-M8} and Proposition \ref{C:prop-M8} the manifold $M^8$ satisfies the assumptions of Proposition \ref{prop:FinProps-inf-family} for $k=4$. Theorem \ref{ssgrofhyperbolicgroup} is an immediate consequence.
\hfill $\Box$


\bigskip
\bigskip
\begin{tiny}
\begin{tabular}{lllllll}
Claudio Llosa Isenrich & & & Bruno Martelli & & & Pierre Py \\
Faculty of Mathematics & & & Dipartimento di Matematica & & & IRMA \\
Karlsruhe Institute of Technology & & & Largo Pontecorvo 5 & & & Universit\'e de Strasbourg \& CNRS \\
76131 Karlsruhe, Germany  & & & 56127 Pisa, Italy & & & 67084 Strasbourg, France \\
claudio.llosa at kit.edu & & & bruno.martelli at unipi.it & & & ppy at math.unistra.fr \\    
\end{tabular}

\end{tiny}

\end{document}